\documentclass[preprint,12pt]{amsart}

\usepackage{amssymb}
\usepackage[english]{babel}
\usepackage{amsmath}
\usepackage{amsfonts}
\usepackage{amssymb}
\usepackage{graphicx}

\newtheorem{Theorem}{Theorem}[section]

\newtheorem{Lemma}[Theorem]{Lemma}
\newcommand{\ds}{\displaystyle}

\usepackage{color}

\usepackage[dvips]{epsfig}
\begin{document}



\title[A locking-free sparse optimal control problem of Timoshenko beam]{A locking-free optimal control problem with $L^1$ cost for optimal placement of control devices in Timoshenko beam} \footnote{Partially supported by Conicyt-ANILLO ACT1106, Conicyt-REDES 140183, MATHAMSUD Project 15MAT-02, FONDECYT 1140392, EPN Project PIS-14-13}  
\author{E. Hern\'andez}
\author{P. Merino}

\address[Erwin Hern\'andez]{Departamento de Matem\'atica, Universidad T\'ecnica Federico Santa Mar\'ia, Casilla 110-V, Valpara\'iso, Chile. }
\address[Pedro Merino]{Centro de Modelaci\'on Matem\'atica (MODEMAT), EPN--Quito, Ladr\'on guevara E11-253, Quito, Ecuador.}

\begin{abstract}

\noindent The numerical approximation of an optimal control problem with $L^1$-control of a Timoshenko beam is considered and analyzed by using the finite element method. From the practical point of view, inclusion of the $L^1$--norm in the cost functional is interesting in the case of beam vibration model, since the sparsity enforced by the $L^1$--norm is very useful for localizing actuators or control devices. The discretization of the control variables is performed by using piecewise constant functions. The states and the adjoint states are approximated by a \emph{locking free scheme} of linear finite elements. Analogously to the purely $L^2$--norm penalized optimal control, it is proved that this approximation have optimal convergence order, which do not depend on the thickness of the beam.
\end{abstract}

\maketitle


\section{Introduction}
The optimal control and stability properties of flexible
beams have been considered by many researchers in the last years,
mainly motivated for science and engineering applications (see, for
instance, \cite{fuller}). A relevant fact in this area is to develop
efficient numerical methods for both: the control problem and the
beam structure model (seen as state equation), respectively.

The numerical analysis of optimal control problems has been an active research area, in particular in the derivation of a priori error estimates arising in its numerical approximation. 
Also, the analysis of problems involving a functional that contains an $L^1(\Omega)$--control cost term
has been considered in the literature. In \cite{Casas2017} we can find a review on sparse control for differential equations. The article \cite{Stadler} seems to be the first to provide an analysis when the distributed control problem associated to a linear elliptic equation is considered. The author utilizes a regularization technique that involves an $L^2(\Omega)$--control cost term and analyze optimality conditions and study the convergence properties of a semismooth Newton method for computing numerically the optimal control. These results were later extended in \cite{Wachsmuth2}, where the authors obtain rates of convergence with respect to a regularization parameter. Subsequently, in \cite{CasasHerzogWach2}, the authors consider a semilinear and elliptic PDE as state equation and analyze second order optimality conditions. Simultaneously, the numerical analysis based on finite element discretizationis considered in \cite{Wachsmuth2}, where the state equation is a linear and elliptic PDE and in \cite{CHW:12again,CasasHerzogWach2} where an extensions to the  semilinear case is introduced.  This kind of problem are also extended by considering the parabolic case in \cite{CasasHerzog2017} and fractional diffusion state equation in \cite{OtarolaSalgadoSparse}.

On the other hand, the most common mathematical model used for thick beam is the
Timoshenko model.  In \cite{Beams} it is concluded that
Timoshenko model is remarkably more accurate if it is compared with
other theories of beam structures (for instance, Euler-Bernoulli
model). Nevertheless, it is well-known that standard finite element
methods applied to this model produce very unsatisfactory results
when the thickness of the beam goes to zero; this fact is known as
locking phenomenon (see \cite{arnold}, \cite{hors}). From the numerical
analysis point of view, this phenomenon can be appreciated in a
priori error estimates for the method considered, because the
associate constants depend on the thickness of the structure in such
a way that they degenerate when this parameter becomes small, this is an important drawback when the solution is considered to be controlled. To
avoid numerical locking, special methods based on reduced
integration or mixed formulations have been considered and mathematically studied. The paper \cite{arnold}, is the first work in which has been proved
that locking arises because of the shear term and has been proposed
and analyzed a locking-free method based on a mixed formulation.
This proposed method has been used and analyzed when it is
applied to the problem of free vibration of a general curved rod
(see \cite{hors}), which covers the Timoshenko beam case. 

The mathematical analysis of the optimal control problems of Timoshenko beams has been considered in \cite{MM}),
despite this, the numerical analysis point of view is considered in \cite{HO1,HO2}, where an active control vibration problem
is studied. Moreover, to the best of the author knowledge, in all references 
the localization of the actuator is fixed and choices without any realistic considerations.

 Although the the analysis of the optimal control problem of Timoshenko beams is covered by the theory of \cite{Stadler} in a general setting, the derivation of the numerical results for non differentiable optimal control problems is not a trivial task. There are only a few contributions from the practical point of view of sparse optimal control problems, therefore it is important revisiting the analysis for specific models including engineering purposes. In particular, in the case of Timoshenko beam model, the incorporation  of $L^1$--norm cost has very interesting potential applications in the placement of control devices. 
 The contribution of this work focuses in deriving convergence results by considering a locking free numerical approximation for the optimization of Timoshenko beam model $L^1$ sparsity inducting term in the cost functional.  The lack of differentiability of this cost term entails new numerical challenges in their approximation and numerical resolution. Therefore, we discuss these topics for this type of problems and propose numerical methods for its numerical solution.

The paper is organized as follows: In Section 2 we present the problem considered in the continuous form, which is fully
discretized in Section 3, where we introduce the locking free finite element methods applied to the beam, and the discretization of the
control problem. Finally, in Section 4 we include numerical examples to show the theoretical results.

 \section{The optimal control problem}

In the following, we shall describe the continuous problem that will
be analyzed and give a brief description on their properties. 
Throughout this paper, $C$ denotes a strictly positive constant, not necessarily the same at each occurrence, but always independent of
the thickness $t$ and of the mesh-size $h$.

Let's consider the following optimal control problem
\begin{equation}
\label{J} \textrm{minimize }J(w,u)=
\frac{1}{2}\|w-w_d\|^2_{L^2(\mathbf{I})}+\frac{\nu}{2}\|u\|^2_{L^2(\mathbf{I})}+\frac{\eta}{2}\|u\|_{L^1(\mathbf{I})}
\end{equation}
subject to the Timoshenko equations (state equations)
\begin{eqnarray}\label{eq1}
\left\{ \begin{array}{rcll}
\ds kAG\left(\frac{d^2w}{dx^2}-\frac{d \theta }{dx}\right)&=&f+u\quad x\in \mathbf{I},\vspace{0.15cm}\\
\ds EI\frac{d^2\theta}{dx^2}+kAG\left(\frac{dw}{dx}-\theta
\right)&=&g\quad x \in \mathbf{I}, \vspace{0.1cm} \\
w(0)=w(L)=\theta(0)=\theta(L)&=&0,\\
 \end{array}
\right.
\end{eqnarray}
and subject to the control constrains
\begin{eqnarray}
\label{cc} u_a\leq u(x) \leq u_b, \quad \textrm{for a.a  } x
\in \mathbf{I},
\end{eqnarray}
where $w$ denote the transversal displacement of the beam, $\theta$
the rotation of its midplane and $\mathbf{I}:=(0,L)$, with $L$ the
length of the beam. The elastic beam of thickness $t\in(0,1]$ has a
reference configuration $\mathbf{I} \times (-t/2,t/2)$. The
coefficients $E$ and $I$, that will be assumed constants, represent
the Young modulus and the inertia moment,
respectively. The coefficient $k$ is a correction factor usually
taken as $5/6$; $A$ and $G$ represent the sectional area of the beam
and elasticity modulus of shear. The term $f$ represents an extern
load and $g$ the bending moment. Moreover, $\nu>0$  and $\eta>0$ represents the
cost of control, $u_a$ and $u_b$ are function in $L^\infty(\mathbf{I})\cap H^1(\mathbf{I})$ and $w_d, f$ and
$g$ are given functions in $L^2(\mathbf{I})$. Note that we consider
two control forces: on the transversal displacement and on the
rotation of the midplane of the beam.

The set of admissible controls is given by $U_{ad}$:
$$
U_{ad}:=\left\{ u \in L^2(\mathbf{I}): u_a \leq u \leq u_b, \textrm{
a.e. } x \in \mathbf{I}\right\}.
$$

We will adapt and extend the convergence results presented in the work of
Wachsmuth and Wachsmuth \cite{Wachsmuth2} to our case. Our aim is obtaining convergence results, independent of the thickness of the beam, for the finite element approximation of the optimal control.

We begin our analysis with the standard fact that for every $f,g \in
L^2(\mathbf{I})$, the unique solution $(w,\theta)$ of the problem
\begin{eqnarray}\label{eqadc}
\left\{ \begin{array}{rcll}
\ds kAG\left(\frac{d^2w}{dx^2}-\frac{d \theta }{dx}\right)&=&f \quad x\in \mathbf{I},\vspace{0.15cm}\\
\ds EI\frac{d^2\theta}{dx^2}+kAG\left(\frac{dw}{dx}-\theta
\right)&=&g \quad x \in \mathbf{I}, \vspace{0.1cm} \\
w(0)=w(L)=\theta(0)=\theta(L)&=&0,\\
 \end{array}
\right.
\end{eqnarray}
belongs to $H_0^1(\mathbf{I})^2 \cap H^2(\mathbf{I})^2$ (see, for
instance, \cite{arnold}). Moreover, there exists a positive constant
$C$ such that
$$
\|(w,\theta)\|_{H^2(\mathbf{I})^2}\leq
C\|(f,g)\|_{L^2(\mathbf{I})^2}.
$$
Now, it is necessary to introduce the adjoint problem
\begin{eqnarray}\label{eqad}
\left\{ \begin{array}{rcll}
\ds kAG\left(\frac{d^2p}{dx^2}-\frac{d q }{dx}\right)&=&w-w_d \quad x\in \mathbf{I},\vspace{0.15cm}\\
\ds EI\frac{d^2q}{dx^2}+kAG\left(\frac{dp}{dx}-q
\right)&=&\theta\quad x \in \mathbf{I}, \vspace{0.1cm} \\
p(0)=p(L)=q(0)=q(L)&=&0,\\
 \end{array}
\right.
\end{eqnarray}
It is clear that the adjoint problem admits a unique solution $(p,q)
\in H_0^1(\mathbf{I})^2 \cap H^2(\mathbf{I})^2$, which is embedded
continuously in $\mathcal{C}^{0,1}(\overline{\mathbf{I}})$. In addition, the existence of a positive constant
$C$ is guaranteed, such that
$$
\|(p,q)\|_{H^2(\mathbf{I})^2}\leq
C\|(w-w_d,\theta)\|_{L^2(\mathbf{I})^2}.
$$

The weak formulations associated to the problems (\ref{eqadc}) and
(\ref{eqad}), respectively, are written in the following manner:

\textit{Find $(w_{t},\theta_{t}) \in H_0^1(\mathbf{I})^2$ such that}
\begin{equation}\label{PV1}
 \ds \frac{E}{12} \int_{\mathbf{I}}\frac{d \theta_{t}}{dx}
\frac{d\beta}{dx} dx+
 \frac{\kappa}{t^2}\int_{\mathbf{I}}\left(\ds \frac{dw_{t}}{dx}- \theta_{t}\right)\left(\frac{dv}{dx}- \beta\right)dx =
\ds\int_{\mathbf{I}} (f+u)v dx + \ds \frac{t^2}{12}
\int_{\mathbf{I}} g \beta dx,\quad \forall(v,\beta) \in H_0^1(\mathbf{I})^2
\end{equation}

and \\

\textit{Find $(p_{t},q_{t}) \in H_0^1(\mathbf{I})^2$ such that}
\begin{equation}\label{PV2}
 \ds \frac{E}{12} \int_{\mathbf{I}}\frac{d q_{t}}{dx}
\frac{d\beta}{dx} dx+ 
\frac{\kappa}{t^2}
 \int_{\mathbf{I}}\left(\ds \frac{dp_{t}}{dx}- q_{t}\right)\left(\frac{dv}{dx}-
\beta\right)dx = \ds\int_{\mathbf{I}} (w_t-w_d)v dx + \ds
\frac{t^2}{12} \int_{\mathbf{I}} \theta_t \beta dx,\quad\forall(v,\beta) \in H_0^1(\mathbf{I})^2
\end{equation}

Also, we consider the bilinear form that appear in the right hand side of the above equations (respectively, changing the variables) :
\begin{equation}
a_t((w_{t},\theta_{t}),(v,\beta) ):=
 \ds \frac{E}{12} \int_{\mathbf{I}}\frac{d \theta_{t}}{dx}
\frac{d\beta}{dx} dx+
 \frac{\kappa}{t^2}\int_{\mathbf{I}}\left(\ds \frac{dw_{t}}{dx}- \theta_{t}\right)\left(\frac{dv}{dx}- \beta\right)dx.
\end{equation}

Hereafter, we denote by $u$ the control function
associated to the problem (\ref{J})-(\ref{cc}), and we call the
solution $(w,\theta)$ of problem (\ref{eq1}) for a given
control $u$, an {\it associated state} to $u$ and write
$(w(u),\theta(u))$. In the same way, we call the solution $(p,q)$,
of the problem (\ref{eqad}) corresponding to $(w(u),\theta(u))$, an
{\it associated adjoint state} to $u$ and write $(p(u),q(u))$. Without loss of generality, we will drop the subindex $t$ and the variable depending $u$ when there is not risk of confussion. 

The corresponding solution mapping is denoted by $\mathcal S$. It is clear that $\mathcal S$ is a continuous linear injective operator from $H^{-1}({\mathbf{I}})$ to $H_0^1({\mathbf{I}})$, such that $w={\mathcal S}(f+u)$. The adjoint operator of $S$ will be denoted by ${\mathcal S}^*$.

The existence and uniqueness of the solutions of our  optimal control problem (\ref{J})-(\ref{cc}) follow from the analysis in Section 2 of \cite{Wachsmuth2}, by considering that ${\mathcal S}$ is injective. The next step in our analysis is formulating the necessary and sufficient optimality conditions for the solution of our problem. These conditions are stablished in the next Lemma. 

\begin{Lemma}\label{l:fonc} Let  $\bar{u}\in U_{ad}$  with associated state $(w(\bar{u}),\theta(\bar{u}))$ the solution of problem (\ref{J})-(\ref{cc}). Then, there exists an adjoint state $p(\bar{u})={\mathcal S}^*(w(\bar{u})-w_d)$, and a subgradient $\lambda$ in the subdiferential of $\eta\|\bar{u}\|_{L^1({\mathbf I})}$ satisfying the following variational inequality
\begin{equation}
\left(-{p}(\bar{u})+\nu \bar{u}+\lambda,
v- \bar{u}\right)\geq 0
\qquad \forall v \in U_{ad}.\label{contIneq}
\end{equation}
\end{Lemma}
\textit{Proof.} See Section 2 in \cite{Wachsmuth2}
\qed

A well known but important observation is that problem (\ref{J})-(\ref{cc}) is a convex optimization problem. Therefore, first order conditions derived in Lemma \ref{l:fonc} are also sufficient.
\section{Discretization and convergence results}

This section is devoted to the analysis of afull discretization of problem (\ref{J})-(\ref{cc}) by using the finite element method. As pointed out, to avoid dependence of the convergence properties on the thickness parameter we consider the locking free scheme proposed in \cite{arnold}. In order to derive the
numerical approximation of the discrete state adjoint
states, we will prove that this result does not present numerical locking. First, we will indicate the reason why the introduction of the modified locking free method is necessary.

\subsection{Fully discretized problem}

The following step is the discretization of the optimal control
problem (\ref{J})-(\ref{cc}).  Let's
consider $\{ \mathcal{T}_h\}$ a family of partitions of the interval
$\mathbf{I}$:
$$
\mathcal{T}_h: 0=s_0<s_1<\cdots s_n=L,
$$
with mesh-size
$$
h:=\max_{j=1,\ldots,n}\left(s_j-s_{j-1} \right).
$$

The control variables $u$ will be discretized by piecewise
constant elements on the mesh $\mathcal{T}_h$ using the following
discrete space:
$$
U_h:=\left\{ u \in L^{\infty}(\mathbf{I}): u|_{[s_{j-1},s_j]} \in
\mathbb{P}_0, j=1,\ldots,n\right\}\subset L^2(\mathbf{I}).
$$

For the beam solution, we consider the following finite element space:
$$
V_h:=\left\{ v \in H_0^1(\mathbf{I}): v|_{[s_{j-1},s_j]} \in
\mathbb{P}_1, j=1,\ldots,n\right\}\subset H_0^1({\mathbf{I}}),
$$where $\mathbb{P}_1$ denotes the space of polynomials of grade less
than or equal to $1$.

With these definitions, the standard procedure is to consider the discrete weak formulation of the solution operator $\mathcal S$, i.e., we consider the following finite dimensional variational problem:

\textit{Find $(w_{th},\theta_{th}) \in V_h^2$ such that}
\begin{align}\label{PVd}
& \ds \frac{E}{12} \int_{\mathbf{I}}\frac{d \theta_{th}}{dx}
\frac{d\beta_h}{dx} dx+
 \frac{\kappa}{t^2}\int_{\mathbf{I}}\left(\ds \frac{dw_{th}}{dx}- \theta_{th}\right)\left(\frac{dv_h}{dx}- \beta_h\right)dx =
\ds\int_{\mathbf{I}} (f+u)v_h dx + \ds \frac{t^2}{12}
\int_{\mathbf{I}} g \beta_h dx, \nonumber\\
& \forall(v_h,\beta) \in V_h^2.
\end{align}
The corresponding discrete solution operator is denoted by ${\mathcal S}_h$. In this context, it is shown in \cite{arnold} that the standard finite elements method applied to the  Timoshenko beam problem (\ref{PV1}) is subjected to the numerical locking phenomenon, this means that they produce unsatisfactory
results for very thin beams. This effect is caused by the shear
stress term. In fact,  if we consider standard finite element
methods for solving (\ref{PV1}), we obtain existence and uniqueness of
the discrete solution $(w_{th},\theta_{th})$ only for  $h<C/t^2$ and
the following (very poor) estimation holds:
\begin{eqnarray*}
\| (w,\theta)-(w_{th},\theta_{th})\|_{L^2({\mathbf{I}})^2} +h\| (w,\theta)-(w_{th},\theta_{th})\|_{H^1({\mathbf{I}})^2}\leq
\frac{C}{t^2}h^2,
\end{eqnarray*}
which implies that we have $$\|{\mathcal S}-{\mathcal S}_h\|_{L^2\longrightarrow H^1}\le \frac{C}{t^2}h.$$

Hence, following the analysis of Section 4 in \cite{Wachsmuth2}, we can prove an a priori finite element error estimation for the approximation of the control function. However, it depends on the thickness of the beam. This effect can be observed in the numerical examples and it represents a serious problem when real control need to be obtained. 

To avoid the numerical-locking in the beam,
Arnold\cite{arnold} introduces and analyzes a locking-free method
based on a mixed formulation of the problem; there,  it was also proved that
this mixed method  is equivalent to using a reduced--order scheme for
the integration of the shear terms in the primal formulation. These
ideas have been extended to the vibration modes of a Timoshenko
curved rod with arbitrary geometry in \cite{hors}, and has been used in optimal control problem for vibration in  \cite{HO1} as well as in the steady case c.f. \cite{HO2}.
\\

In order to apply a mixed locking free scheme to the Timoshenko
equations, we also consider the space
$$
W_h:=\left\{ \frac{dv}{dx} +c :~v \in \mathcal{V}_h, c \in
\mathbb{R}\right\}\subset L^2(\mathbf{I}).
$$

We denote by $\left(w_{h}(u_h), \theta_{h}(u_h)\right)$ the unique element on
$V_h^2$ satisfying the following mixed problem:
\textit{Find $(w_{h},\theta_{h},\gamma_{h}) \in V_h^2 \times W_h$
such that}
\begin{eqnarray}
\label{MFD} \left\{
\begin{array}{rcll}
 \ds \frac{E}{12} \int_{\mathbf{I}}\frac{d \theta_{h}}{dx}
\frac{d\beta_h}{dx} dx+
\int_{\mathbf{I}}\gamma_{h}\left(\frac{dv_h}{dx}- \beta_h\right)dx
&=& \ds\int_{\mathbf{I}} (f+u_{h})v_h dx \\
& & + \ds \frac{t^2}{12} \int_{\mathbf{I}} g\beta_h dx,
\\
\ds \frac{t^2}{\kappa}\int_{\mathbf{I}}\gamma_{th}\eta_h dx&=&\ds
\int_{\mathbf{I}}\left(\ds \frac{dw_{h}}{dx}-
\theta_{h}\right)\eta_h dx,
\end{array}
\right.
\end{eqnarray}
\textit{for all $(v_h,\beta_h) \in V_h^2$ and for all $\eta_h \in
W_h$ respectively.}

Analogouslyº, the adjoint equation is discretized in the same
way and the associate mixed problem is written:

\textit{Find $(p_{h},q_{h},r_{h}) \in V_h^2 \times W_h$
such that}
\begin{eqnarray}
\label{MFDA} \left\{
\begin{array}{rcll}
 \ds \frac{E}{12} \int_{\mathbf{I}}\frac{d q_{h}}{dx}
\frac{d\beta_h}{dx} dx+
\int_{\mathbf{I}}r_{h}\left(\frac{dv_h}{dx}- \beta_h\right)dx &=&
\ds\int_{\mathbf{I}} (w_{h}-w_d)v_h dx \\& & + \ds \frac{t^2}{12}
\int_{\mathbf{I}} (\theta_{h}-\theta_d) \beta_h dx,
\\
\ds \frac{t^2}{\kappa}\int_{\mathbf{I}}r_{h}\eta_h dx&=&\ds
\int_{\mathbf{I}}\left(\ds \frac{dp_{h}}{dx}- q_{h}\right)\eta_h
dx,
\end{array}
\right.
\end{eqnarray}
\textit{for all $(v_h,\beta_h) \in V_h^2$ and for all $\eta_h \in
W_h$ respectively.}

The problem above has been analyzed and error estimates have been
obtained in \cite{arnold} (see also \cite{hors}) and complemented with pointwise error estimates in  
\cite{HO1}. The following results will be used throughout this article, and are summarized in the following
Lemma, which is a direct consequence  of Theorem 3.1 in \cite{hors} and Theorem 4.6 in \cite{HO1}.
\begin{Lemma}\label{lemmafem}
For a given $t>0$, let $(w,\theta)$, $(p,q)$, $(w_{h},\theta_{h})$ and
$(p_{h},q_{h})$ let be the unique solutions of the problems (\ref{PV1}),
(\ref{PV2}), (\ref{MFD}) and (\ref{MFDA}) respectively. Then, the
following estimates hold:
\begin{align}\nonumber
\|(w,\theta)-(w_{h},\theta_{h})\|_{L^2(\mathbf{I})^2} +h\|(w,\theta)-&(w_{h},\theta_{h})\|_{H^1(\mathbf{I})^2} \\
\leq& Ch^2 \bigg( \|u\|_{L^2(\mathbf{I})}+\|f\|_{L^2(\mathbf{I})} +t^2(\|g\|_{L^2(\mathbf{I})}) \bigg),  \label{wt-wth}
\\
\nonumber \|(p,q)-(p_{h},q_{h})\|_{L^2(\mathbf{I})^2} +h \|(p,q)-(p_{h}&,q_{h})\|_{H^1(\mathbf{I})^2} \\
 \leq Ch^2 \bigg(\|u\|_{L^2(\mathbf{I})}+\|f\|_{L^2(\mathbf{I})}&+\|w_d\|_{L^2(\mathbf{I})} +t^2\left(\|\theta_d\|_{L^2(\mathbf{I})} +\|g\|_{L^2(\mathbf{I})}\right)\bigg), \label{pt-pth}\\
\|(w,\theta)-(w_{h},\theta_{h})\|_{L^\infty(\mathbf{I})^2} \leq Ch. \label{feminf}
\end{align}
\end{Lemma}

We will denote by $\hat{\mathcal{S}}_h$ the locking free resolution operator, and the corresponding adjoint operator will be denoted by $\hat{\mathcal{S}}_h^*$. Tehrefore, Lemma \ref{lemmafem} implies that
\begin{equation}\label{S-Sh}\|{\mathcal S}-\hat{\mathcal{S}}_h\|_{L^2\longrightarrow H^1}\le {C}h,\end{equation}
for a positive constant $C$, independent of $t$ and therefore, do not deteriorate when the thickness of the beam goes to zero. 

In order to formulate our discrete version of the optimal control problem, we need to introduce the following  quasi-interpolation operator (see \cite{Wachsmuth2,JC2008} for details). 

Let $\{\phi_i\}_{i=1:N}$ a basis of the discrete space $U_h$, we consider the operator $\Pi_h:L^1(\mathbf{I})\to U_h$, such that 
\begin{equation}\label{Piop}
\Pi_h(u):=\sum_{i=1}^N \left(\frac{1}{s_i-s_{i-1}}\int_{s_{j-1}}^{s_j}u \right)\phi_i,
\end{equation}
which satisfies the following estimate:
\begin{equation}\label{Piop2}
h\|u-\Pi_h(u)\|_{L^2(\mathbf{I})}+\|u-\Pi_h(u)\|_{H^{-1}(\mathbf{I})}\le Ch^2\|\nabla u\|_{L^2(\mathbf{I})^2}.
\end{equation}
Then, we define the discrete admissible set by using the above operator:  
$$
U_{h_{ad}}:=\left\{ u_h \in U_h: \Pi_h(u_a) \leq u \leq \Pi_h(u_b), \textrm{
a.e. } x \in \mathbf{I}\right\}.
$$

With the above definitions, the discrete problem reads

\begin{equation}
\label{Jh} \textrm{minimize }J(u_h)=
\frac{1}{2}\|w_h-w_d\|^2_{L^2(\mathbf{I})}+\frac{\nu}{2}\|u_h\|^2_{L^2(\mathbf{I})}+\frac{\eta}{2}\|u_h\|_{L^1(\mathbf{I})}
\end{equation}
subjet to $u_h\in U_{h_{ad}}$ and (\ref{MFD}). This problem has a unique solution, which is characterized by the following optimality system
\begin{eqnarray}\label{eulercont}\nonumber
w_h&=&\hat{\mathcal{S}}_h(u_h+f)\\\nonumber
p_h&=&\hat{\mathcal{S}}_h^*(w_h-w_d)\\\left(-{p_h}({u_h})+\nu {u_h}+\lambda_h,
v_h- {u_h}\right)&\geq& 0
\qquad \forall v_h \in U_{h_{ad}},\label{discIneq}
\end{eqnarray}
where, $\lambda_h$ is a subgradient in the subdiferential for the $L^1$-norm of $\eta\|u_h\|_{L^1({\mathbf I})}$

Notice that, according to \cite{arnold} (see also \cite{hors,HO1}), the numerical argument to the locking treatment do not involve more computational cost, because this mixed form is equivalent to a reduced integration of the shear term. 

\subsection{Error Estimates}
To derive error estimates we will follow the analysis given in Section 4 of \cite{Wachsmuth2}
First, we note that, in general, ${u}$ does not belong to $U_{h_{ad}}$ and the same is true for  ${u}_h$ in $U_{ad}$.  Therefore, we need to consider ${\tilde u}_h\in U_{h_{ad}}$ as a suitable approximation of $u$, and ${\tilde u}\in U_{ad}$ and approximation of $u_h$, respectively.  

By using ${\tilde u}_h$ and ${\tilde u}$ as test functions, adding the inequalities (\ref{contIneq}) and (\ref{discIneq}), and the definition of the subdiferential,  we have:
\begin{eqnarray*}
\nu \|u-u_h\|^2_{L^2(\mathbf{I})}\le (\nu u_h-p_h,{\tilde u}_h-u)+(\nu u-p,{\tilde u}-u_h)-(p_h-p,u-u_h)\\
+\eta\Big(\|{\tilde u}\|_{L^1(\mathbf{I})}-\|{ u}_h\|_{L^1(\mathbf{I})}+\|{\tilde u}_h\|_{L^1(\mathbf{I})}-\|u\|_{L^1(\mathbf{I})}\Big),
\end{eqnarray*}
and, by a standard argument we arrive to
\begin{eqnarray*}
\nu \|u-u_h\|^2_{L^2(\mathbf{I})}+\|w-w_h\|^2_{L^2(\mathbf{I})}
&\le& (\nu u-p,{\tilde u}-u_h+{\tilde u}_h-u)\\
&&+ \nu(u_h-u,{\tilde u}_h-u)\\
&&-(w_h-w,(\hat{\mathcal{S}}_h-{\mathcal{S}}){\tilde u}_h+{\mathcal{S}}({\tilde u}_h-u))\\
&&-(w-w_d,(\hat{\mathcal{S}}_h-{\mathcal{S}})({\tilde u}_h-u_h))\\
&&+\eta\Big(\|{\tilde u}\|_{L^1(\mathbf{I})}-\|{ u}_h\|_{L^1(\mathbf{I})}+\|{\tilde u}_h\|_{L^1(\mathbf{I})}-\|u\|_{L^1(\mathbf{I})}\Big).
\end{eqnarray*}

Now, by using (\ref{S-Sh}), we obtain
\begin{eqnarray*}
\frac{\nu}{2} \|u-u_h\|^2_{L^2(\mathbf{I})}+\frac12\|w-w_h\|^2_{L^2(\mathbf{I})}
&\le&\|\nu u-p\|_{H^1(\mathbf{I})} (\|{\tilde u}-u_h\|_{H^{-1}(\mathbf{I})}+\|{\tilde u}_h-u\|_{H^{-1}(\mathbf{I})})\\
&&+\nu \|{\tilde u}_h-u\|^2_{L^2(\mathbf{I})}\\
&&+ Ch^4\|{\tilde u}_h\|^2_{L^2(\mathbf{I})}+\|{\mathcal{S}}\|^2_{{\mathcal L}({H^{-1}(\mathbf{I})},{L^2(\mathbf{I})})}\|{\tilde u}_h-u\|^2_{H^{-1}(\mathbf{I})}\\
&&+Ch^2\|w-w_h\|_{L^2(\mathbf{I})}(\|{\tilde u}_h-u\|_{L^{2}(\mathbf{I})}+ \|u-u_h\|_{L^2(\mathbf{I})})\\
&&+\eta\Big(\|{\tilde u}\|_{L^1(\mathbf{I})}-\|{ u}_h\|_{L^1(\mathbf{I})}+\|{\tilde u}_h\|_{L^1(\mathbf{I})}-\|u\|_{L^1(\mathbf{I})}\Big),
\end{eqnarray*}

which depends on the election of ${\tilde u}_h\in U_{h_{ad}}$ and ${\tilde u}\in U_{ad},$.

\begin{Theorem}[Main Result] There  exists a positive constant $C$
independent of $t$ and $h$ such that
\begin{eqnarray*}
\|u-u_h\|_{L^2(\mathbf{I})^2}
\leq C(h\nu^{-1} +(h^2 \nu^{-3/2}).
\end{eqnarray*}
\end{Theorem}
Proof. See Section 4.2 of \cite{Wachsmuth2}. In that case,  ${\tilde u}_h=\Pi_h(u)$ and   ${\tilde u}$ choosing conveniently.\qed

\section{Numerical solution and Examples}
Several different approaches can be used in order to numerically solve Problem \eqref{J} in an efficient way. Following the functional approach from \cite{Stadler}, we describe briefly the application of the Semi--smooth Newton method (SSN), which was used in our experiments. In practice, SSN algorithm works very well although a very large system of equations must be solved. Alternatively, if the problem is first discretized, it can be solved by descend methods that will reduce the size of the system. In particular, a second order method from \cite{dlrlm17} is very efficient for solving optimal control problems.

 Based on Lemma \ref{l:fonc}, is a standard procedure writting the optimality system as follows. Let $\bar u \in L^2(\mathbf I)$ the optimal control for Problem \eqref{J} with associated optimal state $(\bar w,\bar \theta) \in H_0^1(\mathbf I) \times  H^2(\mathbf I)$, there exist an adjoint state $(p,q) \in H_0^1(\mathbf I)^2 \cap  H^2(\mathbf I)^2$ and multipliers $\lambda, \lambda_a, \lambda_b \in L^2(\mathbf I)$ such that these quantities satisfy the following optimality system.
 
 \begin{align}
 	&\bar w = \mathcal{S} \bar u + \mathcal S f\\
 	& p  = \mathcal{S^*} (\bar w - w_d)\\
 	& -p + \nu \bar u + \lambda + \lambda_b - \lambda_a= 0 \\
 	&  \lambda_b \geq 0, \quad b-\bar u \geq 0, \quad  \lambda_b(b -\bar u) =0 \\
 	&  \lambda_a \geq 0, \quad \bar u -a \geq 0, \quad  \lambda_a(\bar u -a) =0 \\
 	& \lambda(x) = \eta \, \text{sign}(\bar u (x)), \quad  \text{ where } \{x: \bar u (x) \not =0 \}\\
 	& |\lambda(x)| \leq \eta,\quad \text{ where } \{x: \bar u (x) =0 \}
 \end{align}
 
 which can be rewritten more compactly using the $\emph{max}$ and $\emph{min}$ functions and setting $ \mu = \lambda +\lambda_b - \lambda_a$, giving the system
 
 \begin{align}
 	&\bar w = \mathcal{S} \bar u + \mathcal S f, \\
 	& p  = \mathcal{S^*} (\bar w - w_d),\\
 	&  -p + \nu \bar u + \mu = 0 ,\\
 	& \mathcal{C}(\bar u,\mu) = 0,
   \end{align}

where 
\begin{align}
\mathcal{C}(\bar u,\mu):=& \nu\bar u - \max (0, \nu\bar u + \mu -\eta)  - \min (0, \nu\bar u + \mu+ \eta) \nonumber \\
&+ \max (0, \nu(\bar u -b) + \mu-\eta) + \min (0, \nu(\bar u -a) + \mu -\eta)	\nonumber
 \end{align}
 collects all information from the multipliers. From this optimality system, we obtain the following Newton system

\begin{align}
\left[
\begin{array}{llll}
 (I\,\, 0) & -\mathcal S & 0 & 0 \\	
0 & \nu I & (-I \,\, 0) & I \\
(-\mathcal{S^*} \,\, 0)& 0 & (I\,\, 0) & 0 \\
0 & I - \chi_A & 0 & \nu^{-1} \chi_A
\end{array}
\right]
&
\left[
\begin{array}{c}
\left(
\begin{matrix} \delta_w \\ \delta_\theta \end{matrix}
\right)\\
\delta_u \\
\left(
\begin{matrix} \delta_p \\ \delta_q \end{matrix}
\right)\\
\delta_{\mu}
\end{array}
\right]=-\left[ 
\begin{array}{c}
 w  -\mathcal{S}  u 
\\
-p + \nu  u + \mu \\
p - \mathcal{S}^*(w-w_d) \\
\mathcal{C}(u,\mu)
\end{array}
\right],
\end{align}

where $A$ corresponds to the active set, given by $$A= \{ x \in \mathbf I: \nu a < p + \eta \leq 0\} \cup \{ x \in \mathbf I: 0 \leq p - \beta < \nu b \}$$

We report the results of several numerical tests that illustrate different scenarios. The optimization problem was solved by applying the SSN algorithm with the implementation of the locking-free finite element scheme described above, coded in \textsc{matlab}. We used a reduced-order scheme for
the integration of the shear term in the primal formulation, such as
the scheme proposed in \cite{arnold}. As mentioned, this approach is equivalent to the
mixed formulation. 

The physical parameters and the control parameters used in the
numerical resolution of the tests problems are the following:


\begin{center}
\begin{tabular}{l|l}
Elastic moduli: $E=$1.44 $\times 10^{9}$Pa &  Poisson coefficient: $\bar{\nu}=$0.35, \\
Correction factor: $k=5/6$ & Density: $\rho=$7.7 $\times 10^{3}$Kg/$\textrm{m}^3$, \\
\end{tabular}
\end{center}


\subsection{Comparison of $L^2$--controls with sparse controls}

This example is intended to show how sparse controls act in a``located" fashion with respect to $L^2$--controls, which are distributed in the whole domain. In the next example, the following parameters for optimization where chosen. 

\begin{center}\label{f:comp1}
\begin{tabular}{l|l}
 $\nu=5\times 10^{-9}$, &$\eta \in [0, 2.7\times 10^{-5}] $ \\
  $a=-11.05$ & $b=11.05$,\\
 $f(x)=100\sin(8\pi\,x)$ &   $g=0$.
\end{tabular}
\end{center}

We run this example in mesh with 601 nodes and compare the solutions for different values of $\eta:=0:2\times 10^{-5}$.  We summarize the results in Table \ref{t:e1} where we show the different cost and the corresponding $L^2$--norm for each solution, up to the value of $\eta$ where the optimal control becomes zero. For example, the first two solutions are depicted in Figure  \ref{f:comp1}. These plots clearly show the effect of the $L^1$ penalization term. Indeed, in contrast with the pure $L^2$--control ($\eta=0$), we observe that the sparse optimal control ( $\eta=3 \times 10^{-6}$) is nonzero in two well identified sections of the beam, where both controls are active. On the other hand, we can see that the corresponding states, representing the vertical displacements are both close to 0. Although the sparse controlled state has a slightly higher amplitude, the optimal control which produces it has a much lesser $L^2$--norm cost. 

\begin{table}[hbt]\label{t:e1}
\centering
  \begin{tabular}{|l|l|l|l|}
\hline
$\eta$& Cost & $L^2$--norm & Null \\ \hline \hline
0&1.6986e-06&9.4704&0\\\hline
3e-06&6.3031e-06&3.179&530\\\hline
6e-06&8.9758e-06&2.813&545\\\hline
9e-06&1.1125e-05&2.5228&555\\\hline
1.2e-05&1.2841e-05&2.203&564\\\hline
1.5e-05&1.4146e-05&1.8141&571\\\hline
1.8e-05&1.5049e-05&1.2875&576\\\hline
2.1e-05&1.5553e-05&0.66013&582\\\hline
2.4e-05&1.5674e-05&0.046107&596\\\hline
2.7e-05&1.5677e-05&0&600\\\hline
\end{tabular}

  \caption{Solution costs for different values of $\eta$ }
\end{table}

\begin{figure}[h]\label{f:op_ex1}
\centering
  \includegraphics[scale=0.35]{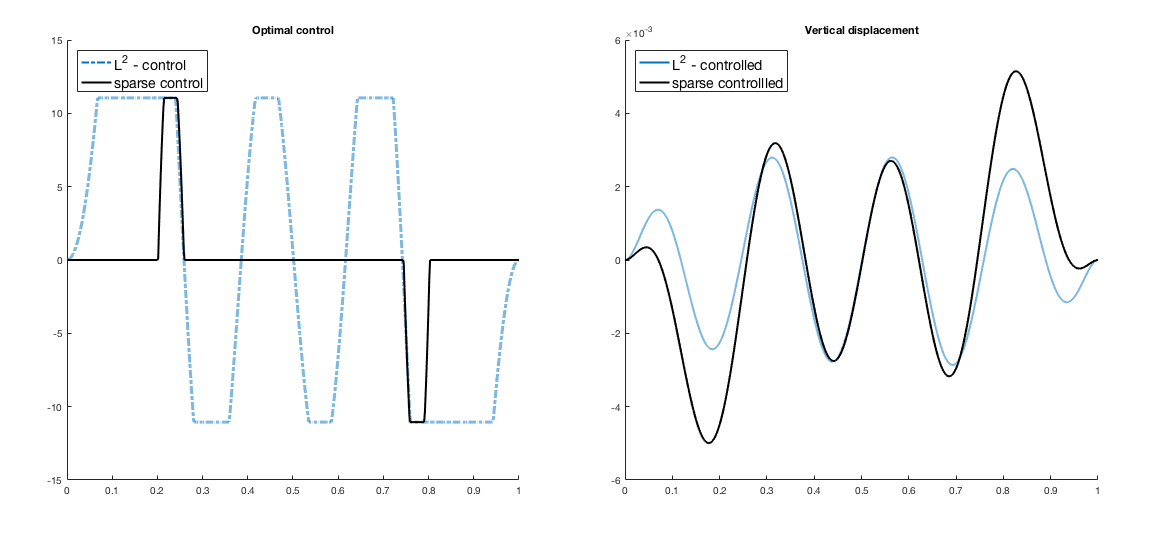}
  \caption{$L^2$--norm control vs sparse control and corresponding optimal states}
\end{figure}

\subsubsection{Testing the locking--free property}

\begin{figure}[hbt]\label{f:lockingfree}
  \includegraphics[scale=0.28]{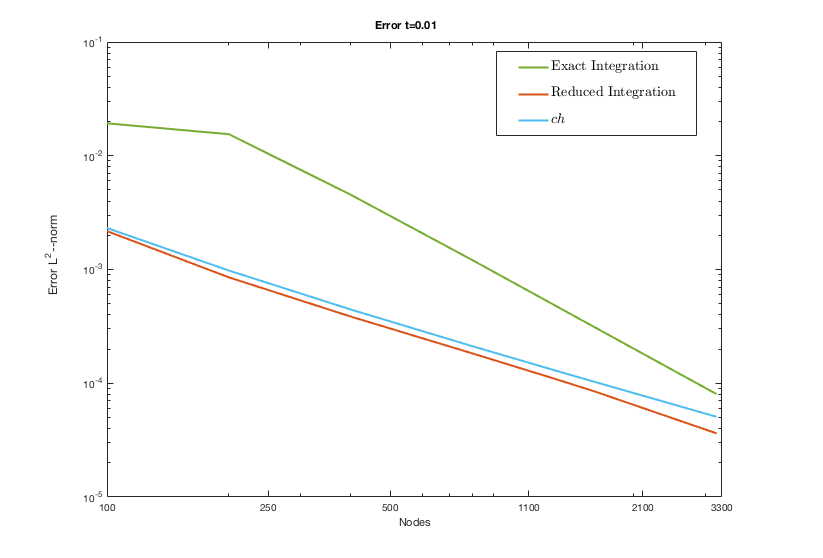}
  \includegraphics[scale=0.28]{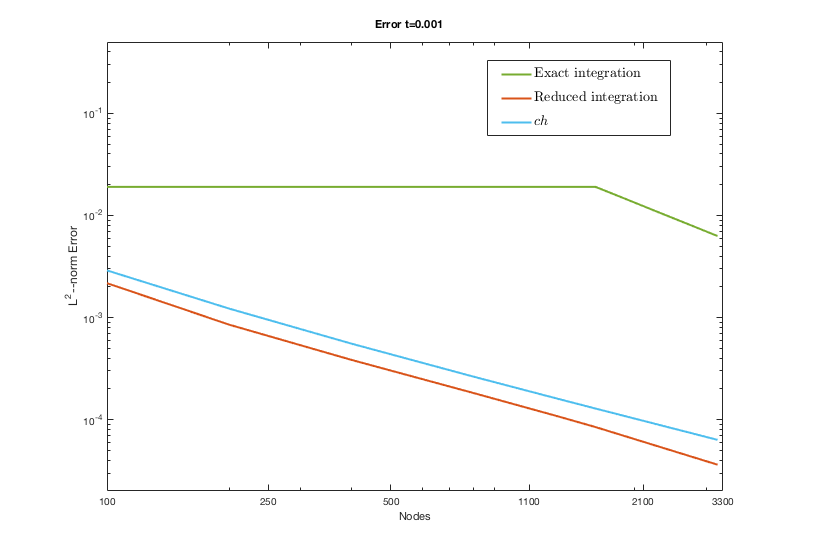}
  \caption{$L^2$--norm error of for optimal control for different thickness parameter}
\end{figure}

In order to illustrate the locking--free feature of our scheme, we consider two
small values of the thickness: $t=10^{-2}$ and
$t=10^{-3}$ (see \cite{hors}, \cite{HO1}). Figure \ref{f:lockingfree} shows error curves in terms of the refinement parameter $N=1/h$ for controls obtained by means
of the classic method, i.e., as a solution of the problem
(\ref{J}): $(u_1,u_2)$ and for the controls obtained with
the superconvergence step: $(\tilde{u}_1,\tilde{u}_2)$. In this
figure it can be clearly seen that the order of convergence is $O(h)$. For $t=0.01$, we observe that a large number of nodes is needed for the usual scheme in order to achieve the precision of the solution computed with the reduced scheme. In addition, when $t=001$  we observe the locking effect in the usual method versus the locking-free property of our scheme.

%
%
%
%
\bibliographystyle{plain}

\end{document}